\documentclass[12pt]{article}

\usepackage{amsmath}
\usepackage{amsfonts}
\usepackage{amssymb}

\newtheorem{theorem}{Theorem}
\newtheorem{corollary}{Corollary}

\newtheorem{lemma}{Lemma}

\newtheorem{remark}{Remark}

\def\r{\mathbb{R}}

\newenvironment{proof}[1][Proof]{\noindent\textbf{#1.} }{\ \rule{0.5em}{0.5em}}

\begin{document}

\title{\bf Hardy's theorem for the $q$-Bessel Fourier transform}
\date{ }
\author{Lazhar Dhaouadi \thanks{%
Institut Pr\'eparatoire aux Etudes d'Ing\'enieur de Bizerte
(Universit\'e du 7 novembre à Carthage). Route Menzel Abderrahmene
Bizerte, 7021 Zarzouna, Tunisia. \quad\quad\quad\quad\quad\quad
E-mail lazhardhaouadi@yahoo.fr}} \maketitle

\begin{abstract}
In this paper we give a q-analogue of the Hardy's theorem for the
$q$-Bessel Fourier transform. The celebrated theorem asserts that if
a function $f$ and its Fourier transform $\widehat{f}$ satisfying
$|f(x)|\leq c.e^{-\frac{1}{2}
x^2}$ and $|\widehat{f}(x)|\leq c.e^{-\frac{1}{2} x^2}$ for all $x\in\mathbb{%
R}$ then $f(x)=\text{const}.e^{-\frac{1}{2} x^2}$.
\end{abstract}

\section{The q-Foureir Bessel transform}

Throughout this paper we consider $0<q<1$ and we adopt the standard
conventional notations of [2]. We put
\begin{equation*}
\mathbb{R}_q^+=\{q^n,\quad n\in\mathbb{Z}\},
\end{equation*}
and for complex $a$
\begin{equation*}
(a;q)_0=1,\quad (a;q)_n=\prod_{i=0}^{n-1}(1-aq^{i}), \quad
n=1...\infty.
\end{equation*}
Jackson's $q$-integral (see [3]) in the interval $[0,\infty[$ is
defined by
\begin{equation*}
\int_0^\infty f(x)d_qx=(1-q)\sum_{n=-\infty}^\infty q^nf(q^n).
\end{equation*}
We introduce the following functional spaces $\mathcal{L}_{q,1,\nu}$
of even functions $f$ defined on $\mathbb{R}_q^+$ such that
\begin{equation*}
\|f\|_{q,1,v}=\left[\int_0^{\infty}|f(x)|x^{2v+1}d_qx\right]<\infty.
\end{equation*}
The normalized Hahn-Exton q-Bessel function of order $\nu>-1$ (see
[5])is defined by
\begin{equation*}
j_{\nu }(z,q)=\sum_{n=0}^{\infty }(-1)^{n}\frac{q^{\frac{n(n-1)}{2}}}{%
(q,q)_{n}(q^{\nu +1},q)_{n}}z^{n}.
\end{equation*}
It's an entire analytic function in $z$.

\begin{lemma}
For every $p\in\mathbb{N}$, there exist $\sigma_p>0$ for which
\begin{equation*}
|z|^{2p}|j_{\nu}(z,q^2)|<\sigma_p e^{|z|},\quad\forall
z\in\mathbb{C}.
\end{equation*}
\end{lemma}

\begin{proof}
In fact
\begin{align*}
|z|^{2p}|j_{\nu}(z,q^2)|&\leq\frac{1}{(q^2,q^2)_{\infty}
(q^{2\nu+2},q^2)_{\infty}}\sum_{n=0}^\infty q^{n(n-1)}|z|^{2n+2p} \\
&\leq\frac{q^{p(p+1)}}{(q^2,q^2)_{\infty} (q^{2\nu+2},q^2)_{\infty}}%
\sum_{n=p}^\infty q^{n(n-2p-1)}|z|^{2n}.
\end{align*}
Now using the Stirling's formula
\begin{equation*}
n!\sim\sqrt{2\pi n}\frac{n^n}{e^n},
\end{equation*}
we see that there exist an entire $n_0\geq p$ such that
\begin{equation*}
q^{n(n-2p-1)}<\frac{1}{(2n)!},\quad\forall n\geq n_0,
\end{equation*}
which implies
\begin{equation*}
\sum_{n=n_0}^\infty q^{n(n-2p-1)}|z|^{2n}<\sum_{n=n_0}^\infty \frac{1}{(2n)!}%
|z|^{2n}<e^{|z|}.
\end{equation*}
Finally there exist $\sigma_p>0$ such that
\begin{equation*}
\frac{|z|^{2p}|j_{\nu}(z,q^2)|}{e^{|z|}}<\sigma_p,\quad\forall z\in\mathbb{C}%
.
\end{equation*}
This finish the proof.
\end{proof}

\bigskip

The q-Bessel Fourier transform $\mathcal{F}_{q,\nu }$ introduced in
[1,4] as follow
\begin{equation*}
\mathcal{F}_{q,\nu }f(x)=c_{q,\nu }\int_{0}^{\infty }f(t)j_{\nu
}(xt,q^{2})t^{2\nu +1}d_{q}t.
\end{equation*}%
where \ \
\begin{equation*}
c_{q,\nu }=\frac{1}{1-q}\frac{(q^{2\nu +2},q^{2})_{\infty }}{%
(q^{2},q^{2})_{\infty }}.
\end{equation*}%
The following theorem was proved in [1]

\begin{theorem}
Given $f\in\mathcal{L}_{q,1,\nu}$ then we have
\begin{equation*}
\mathcal{F}_{q,\nu}^2(f)(x)=f(x),\quad\forall x\in\mathbb{R}_q^+.
\end{equation*}
\end{theorem}

\begin{proof}
See [1] p 3.
\end{proof}

\section{Hardy's theorem}

The following Lemma from complex analysis is crucial for the proof
of our main theorem.

\begin{lemma}
Let $h$ be an entire function on $\mathbb{C}$ such that
\begin{equation*}
|h(z)|\leq C e^{a|z|^2},\quad z\in\mathbb{C},
\end{equation*}

\begin{equation*}
|h(x)|\leq Ce^{-ax^{2}},\quad x\in \mathbb{R},
\end{equation*}
for some positive constants $a$ and $C$. Then $h(z)=$
Const$.e^{-ax^{2}}$.
\end{lemma}

\begin{proof}
See [6] p 4.
\end{proof}

\bigskip Now we are in a position to state and prove the q-analogue of the
Hardy's theorem

\begin{theorem}
Suppose $f\in\mathcal{L}_{q,1,\nu}$ satisfying the following
estimates
\begin{equation*}
|f(x)|\leq C e^{-\frac{1}{2} x^2},\quad\forall x\in\mathbb{R}_q^+,
\end{equation*}

\begin{equation*}
|\mathcal{F}_{q,\nu} f(x)|\leq C e^{-\frac{1}{2} x^2},\quad\forall x\in%
\mathbb{R},
\end{equation*}
where $C$ is a positive constant. Then there exist $A\in\r$ such
that
\begin{equation*}
f(z)=A.c_{q,\nu}\mathcal{F}_{q,\nu}\left(e^{-\frac{1}{2} x^2}\right)(z),%
\quad\forall z\in\mathbb{C}.
\end{equation*}
\end{theorem}

\begin{proof}
We claim that $\mathcal{F}_{q,\nu} f$ is an analytic function and
there exist $C^{\prime}>0$ such that

\begin{equation*}
|\mathcal{F}_{q,\nu }f(z)|\leq C^{\prime
}e^{\frac{1}{2}|z|^{2}},\quad \forall z\in \mathbb{C}.
\end{equation*}
We have
\begin{equation*}
|\mathcal{F}_{q,\nu }f(z)|\leq c_{q,\nu }\int_{0}^{\infty
}|f(x)||j_{\nu }(zx,q^{2})|x^{2\nu +1}d_{q}x.
\end{equation*}
From the Lemma 1, if $\left\vert z\right\vert >1$ then there exist
$\sigma _{1}>0$ such that
\begin{equation*}
x^{2\nu +1}|j_{\nu }(zx,q^{2})|=\frac{1}{\left\vert z\right\vert ^{2\nu +1}}%
(\left\vert z\right\vert x)^{2\nu +1}|j_{\nu }(zx,q^{2})|<\frac{\sigma _{1}}{%
1+\left\vert z\right\vert ^{2}x^{2}}e^{x|z|},\quad \forall x\in \mathbb{R}%
_{q}^{+}.
\end{equation*}
Then we obtain
\begin{equation*}
|\mathcal{F}_{q,\nu }f(z)|\leq C\sigma _{1}c_{q,\nu }\left[
\int_{0}^{\infty
}\frac{e^{-\frac{1}{2}(x-|z|)^{2}}}{1+\left\vert z\right\vert ^{2}x^{2}}%
d_{q}x\right] e^{\frac{1}{2}|z|^{2}}<C\sigma _{1}c_{q,\nu }\left[
\int_{0}^{\infty }\frac{1}{1+x^{2}}d_{q}x\right]
e^{\frac{1}{2}|z|^{2}}.
\end{equation*}
Now, if $\left\vert z\right\vert \leq 1$ then there exist $\sigma
_{2}>0$
such that%
\begin{equation*}
x^{2\nu +1}|j_{\nu }(zx,q^{2})|\leq \sigma _{2}e^{x},\quad \forall
x\in \mathbb{R}_{q}^{+}.
\end{equation*}
Therefore
\begin{equation*}
|\mathcal{F}_{q,\nu }f(z)|\leq C\sigma _{2}c_{q,\nu }\left[
\int_{0}^{\infty }e^{-\frac{1}{2}x^{2}+x}d_{q}x\right] \leq C\sigma
_{2}c_{q,\nu }\left[
\int_{0}^{\infty }e^{-\frac{1}{2}x^{2}+x}d_{q}x\right] e^{\frac{1}{2}%
|z|^{2}}.
\end{equation*}
Which leads to the estimate (1). Using Lemma 2, we obtain
\begin{equation*}
\mathcal{F}_{q,\nu }f(z)=\text{const}.e^{-\frac{1}{2}z^{2}},\quad
\forall z\in \mathbb{C},
\end{equation*}
and by theorem 1, we conclude that
\begin{equation*}
f(z)=\text{const}.\mathcal{F}_{q,\nu }\left(
e^{-\frac{1}{2}t^{2}}\right) (z),\quad \forall z\in \mathbb{C}.
\end{equation*}
This finish the proof.
\end{proof}

\begin{corollary}
Suppose $f\in\mathcal{L}_{q,1,\nu}$ satisfying the following
estimates
\begin{equation*}
|f(x)|\leq C e^{-p x^2},\quad\forall x\in\mathbb{R}_q^+,
\end{equation*}

\begin{equation*}
|\mathcal{F}_{q,\nu} f(x)|\leq C e^{-\sigma x^2},\quad\forall x\in\mathbb{R}%
,
\end{equation*}
where $C,p,\sigma$ are a positive constant and $p\sigma=\frac{1}{4}$
. We suppose that there exist $a\in\mathbb{R}_q^+$ such that
$a^2p=\frac{1}{2}$. Then there exist $A\in\r$ such that
\begin{equation*}
f(z)=A.c_{q,\nu}\mathcal{F}_{q,\nu}\left(e^{-\sigma
t^2}\right)(z),\quad\forall z\in\mathbb{C}.
\end{equation*}
\end{corollary}

\begin{proof}
Let $a\in\mathbb{R}_q^+$, and put
\begin{equation*}
f_a(x)=f(ax),
\end{equation*}
then
\begin{equation*}
\mathcal{F}_{q,\nu} f_a(x)=\frac{1}{a^{2\nu+2}}\mathcal{F}_{q,\nu}
f(x/a).
\end{equation*}
In the end, applying Theorem 2 to the function $f_a$.
\end{proof}

\bigskip

\begin{remark}
Using the $q$-Central limit Theorem (see [1]) we give a probability
interpretation of the function
$c_{q,\nu}\mathcal{F}_{q,\nu}\left(e^{-\sigma t^2}\right)$ . In fact
if $(\xi_n)_{n\geq 0}$ be a sequence of positive probability
measures of $\mathbb{R}_q^+$, satisfying

\begin{equation*}
\lim_{n\rightarrow \infty}n\sigma_n=\frac{(q^2,q^2)_1(q^{2v+2},q^2)_1}{q^2}%
\sigma,\quad\text{where}\quad
\sigma_n=\int_0^{\infty}t^2t^{2v+1}d_q\xi_n(t),
\end{equation*}
and
\begin{equation*}
\lim_{n\rightarrow \infty}n\widetilde{\sigma}_n=0,
\quad\text{where}\quad
\widetilde{\sigma}_n=\int_0^{\infty}\frac{t^4}{1+t^2}t^{2v+1}d_q\xi_n(t),
\end{equation*}
then the nth $q$-convolution product $\xi_n^{*n}$ converge strongly
toward a measure $\xi$ defined by
\begin{equation*}
d_q\xi(x)=c_{q,\nu}\mathcal{F}_{q,\nu}\left(e^{-\sigma
t^2}\right)(x) d_qx.
\end{equation*}
\end{remark}

\begin{corollary}
Suppose $f\in\mathcal{L}_{q,1,\nu}$ satisfying the following
estimates
\begin{equation*}
|f(x)|\leq C e^{-p x^2},\quad\forall x\in\mathbb{R}_q^+,
\end{equation*}

\begin{equation*}
|\mathcal{F}_{q,\nu} f(x)|\leq C e^{-\sigma x^2},\quad\forall x\in\mathbb{R}%
,
\end{equation*}
where $C,p,\sigma$ are a positive constant and $p\sigma>\frac{1}{4}$
.We suppose that there exist $a\in\mathbb{R}_q^+$ such that
$a^2p=\frac{1}{2}$. Then $f\equiv 0$.
\end{corollary}

\begin{proof}
In fact there exist $\sigma^{\prime}<\sigma$ such that $p\sigma^{\prime}=%
\frac{1}{4}$. Then the function $f$ satisfying the estimates of
Corollary 1, if we replacing $\sigma$ by $\sigma^{\prime}$. Which
implies
\begin{equation*}
\mathcal{F}_{q,\nu}
f(x)=\text{const}.e^{-\sigma^{\prime}x^2},\quad\forall
x\in\mathbb{R}.
\end{equation*}
On the other hand, $f$ satisfying the estimates of Corollary 2, then
\begin{equation*}
|\text{const}.e^{-\sigma^{\prime}x^2}|\leq C e^{-\sigma
x^2},\quad\forall x\in\mathbb{R}.
\end{equation*}
This implies $\mathcal{F}_{q,\nu} f\equiv 0$, and by Theorem 1 we
conclude that $f\equiv 0$
\end{proof}

\begin{remark}
Hardy's theorem asserts the impossibility of a function and its
q-Fourier Bessel transform to be simultaneously "very rapidly
decreasing". Hardy's theorem can also be viewed as a sort of
"Qualitative uncertainty principles". One such example can be the
fact that a function and its $q$-Bessel Fourier transform cannot
both have compact support.
\end{remark}

\end{document}